\documentclass[11pt,a4paper]{article}
\usepackage{amssymb,amscd}
\title{$O(n)$ invariant solutions of Abreu's equation}  

\author{
A.N.W. Hone\thanks{Institute of Mathematics \&
Statistics,
University of Kent, Canterbury CT2 7NF, UK.
E-mail:anwh@ukc.ac.uk}
}

\begin{document}
\renewcommand{\theequation}{\arabic{section}.\arabic{equation}}
\newcommand{\beq}{\begin{equation}}
\newcommand{\eeq}{\end{equation}}
\newcommand{\bea}{\begin{eqnarray}}
\newcommand{\eea}{\end{eqnarray}}
\maketitle

\begin{abstract}
We consider a fourth order partial differential equation 
in $n$ dimensional space 
introduced by Abreu in the context of K\a"{a}hler metrics 
on toric orbifolds. Similarity solutions depending only 
on the radial coordinate in ${\mathbb R}^n$ are determined 
in terms of a second order ordinary differential equation. 
A local asymptotic analysis of solutions in the neighbourhood of  
singular points is carried out. The integrability (or otherwise) 
of Abreu's equation is discussed.  
\end{abstract} 

\section{Introduction} 

In recent work Abreu has considered toric K\a"{a}hler metrics 
on toric varieties \cite{abreu1} or toric orbifolds   
\cite{abreu2}  
of dimension $2n$. Following a construction due to 
Guillemin \cite{guillemin}, each such variety or orbifold is 
completely determined by its moment polytope $\Delta$ 
in ${\mathbb R}^n$, and the scalar curvature $S$ of 
the K\a"{a}hler metric is given by a formula  
\beq 
S=-\partial_j\partial_k g^{jk}  
\label{eq:curv} 
\eeq 
where the matrix elements $g^{jk}$ are functions of 
${\bf x}=(x_1,\ldots,x_n)^{T}\in{\mathbb R}^n$, 
and $\partial_j$ denotes the partial derivative 
$\partial /\partial x_j$ (the summation convention 
is assumed). 
More precisely, in terms of a potential function 
$g({\bf x})=g(x_1,\ldots,x_n)$ 
on ${\mathbb R}^n$, the metric on the interior of the polytope 
$\Delta$ is  
$$ 
ds^2=g_{jk}\,dx_jdx_k,  
$$  
where $g_{jk}$ are the elements of the Hessian matrix ${\bf G}$, i.e. 
$$ 
{\bf G}=(g_{jk}), \qquad g_{jk}=\partial_j\partial_k g, 
$$ 
and in (\ref{eq:curv}) the curvature  
is  determined  in terms of the inverse matrix 
$$ 
{\bf G}^{-1}=(g^{jk}). 
$$ 

In \cite{abreu1} it is shown that the condition for 
extremal toric metrics is that the curvature is an affine 
function of ${\bf x}$, in other words 
\beq 
\partial_j S=constant, \qquad j=1,\ldots,n. 
\label{extrem} 
\eeq 
Abreu constructs such metrics from potentials of the 
form 
\beq 
g=\frac{1}{2}\sum_{l=1}^{d}\ell_l({\bf x})\log\ell_l({\bf x}), 
\label{eq:polys} 
\eeq 
where the $\ell_l({\bf x})$ are 
affine functions which determine the 
facets of the polytope by the equations $\ell_l({\bf x})=0,$ 
$l=1,\ldots,d$. For the simplest case of Guillemin's 
construction \cite{guillemin}, the polytope  is a cuboid in 
${\mathbb R}^n$, so that $d=2n$, and the toric variety is just 
$(S^2)^{\times n}$, i.e.  
$n$ copies of the 2-sphere obtained by attaching an $n$-torus 
$T^n$ to each point in the interior of the cuboid.    

The purpose of this note is to construct 
other types of solution to the equation (\ref{eq:curv}) for 
the case of constant 
scalar curvature, namely    
\beq 
\partial_j\partial_k g^{jk}=-\kappa, \qquad constant. 
\label{eq:abreueq} 
\eeq 
Henceforth we shall follow \cite{donaldson} and 
refer to (\ref{eq:abreueq}) as Abreu's equation, and show that 
for $n>1$ it admits $O(n)$ invariant solutions of a different 
form compared with (\ref{eq:polys}). 

Rather than being geometric, our motivation for considering 
the partial differential equation (PDE) (\ref{eq:abreueq}) 
comes from the theory of integrable systems. Taking 
derivatives of (\ref{eq:polys}) we have 
\beq 
g_{jk}=\frac{1}{2}\sum_{l=1}^d 
\frac{c_{l,j} c_{l,k}}{\ell_l({\bf x})}, 
\label{absum} 
\eeq 
where the constants $c_{l,j}$  
are the coefficients of the affine functions 
$\ell_l({\bf x})$. The expression (\ref{absum}) is 
reminiscent of the potential for 
vanishing rational solutions \cite{krichrat, shiota}  
of the integrable Kadomtsev-Petviashvili (KP) equation, 
\beq 
\frac{\partial}{\partial x}\left( 
4\frac{\partial u}{\partial t_3}-12u\,\frac{\partial u}{\partial x} 
-\frac{\partial^3 u}{\partial x^3}\right) 
-3\frac{\partial^2 u}{\partial t_2^2}=0. 
\label{kp} 
\eeq 
which takes the form   
\beq 
w=-\sum_{l=1}^d\frac{1}{x-x_l(\underline{t})}, \qquad 
u=-2\frac{\partial w}{\partial x},      
\label{kprat} 
\eeq 
where $\underline{t}$ denotes the dependence of the poles 
on the times $t_2,t_3$ (or an infinite sequence of 
such times in the full KP hierarchy). The dynamics of 
the $d$ poles $x_l(\underline{t})$ is governed by the  
integrable   
Calogero-Moser system for $d$ particles. The potential 
$w$ in (\ref{kprat}) is a central object in the Sato 
formulation of KP theory \cite{sato}, since the pseudo-differential 
Lax operator $L$ is constructed with a dressing operator $W$ 
such that  
\beq 
L=W\partial_x W^{-1}, \qquad W=1+w\partial_x^{-1}+\ldots. 
\label{kplax} 
\eeq 
    
The similarity between (\ref{absum}) and the rational KP 
potential (\ref{kprat}) raises the question of whether 
Abreu's equation (\ref{eq:abreueq}) might be integrable in some 
sense. Following the Ablowitz-Ramani-Segur conjecture \cite{ars} that 
all reductions of integrable PDEs should have the 
Painlev\'{e} property, our natural instinct is to 
seek similarity reductions of (\ref{eq:abreueq}) and 
examine the structure of their singularities in the complex plane.    
To avoid begging the question, 
by an integrable PDE we mean one having a Lax pair (ensuring 
solvability by the inverse scattering transform \cite{abcla}) 
and/or infinitely many symmetries \cite{shabat}.

\section{Similarity solutions} 

\setcounter{equation}{0}

For each $n$, Abreu's equation (\ref{eq:abreueq}) has 
$O(n)$ invariant solutions with the potential $g$ being a 
function of the radial distance only, that is 
$$ 
g=g(r), \qquad r=|{\bf x}|. 
$$ 
In that case the Hessian matrix and its inverse are given by 
\beq 
{\bf G}=\frac{f}{r}\,{\bf 1}+\frac{1}{r}\left( 
\frac{f}{r}\right)'{\bf x x}^{T}, 
\qquad 
{\bf G}^{-1}=\frac{r}{f}\,{\bf 1}+\left(  
\frac{1}{r^2f'}-\frac{1}{rf}\right){\bf x x}^{T}, 
\label{eq:radmat} 
\eeq 
where  
$$ 
f=g', 
$$ 
and the prime denotes $d/dr$.       

Substituting the form (\ref{eq:radmat}) of the 
inverse Hessian matrix into (\ref{eq:abreueq}) yields the 
following third order ordinary differential equation (ODE) for 
$f$: 
\beq 
\left(n+r\frac{d}{dr}\right)\left( 
\frac{A'}{r}+rB'+(n+1)B\right)=-\kappa; \qquad 
A=\frac{r}{f}, \, B=\frac{1}{r^2f'}-\frac{1}{rf}. 
\label{eq:sim} 
\eeq 
After an integration, this yields 
\beq 
\frac{A'}{r}+rB'+(n+1)B=\lambda r^{-n}-\frac{\kappa}{n}, 
\label{eq:int} 
\eeq 
for constant $\lambda$. As a second order 
ODE for $f=f(r)$, the equation (\ref{eq:int}) takes  
the explicit form
\beq 
f''=\left(\frac{\kappa}{n}r-\frac{\lambda}{r^{n-1}}- 
\frac{(n-1)}{f}\right)(f')^2+\frac{(n-1)}{r}f'. \label{eq:ode} 
\eeq 
For the purposes of asymptotic analysis, (\ref{eq:ode}) 
may be conveniently rewritten as  
\beq 
\frac{d}{dr}\log\left(\frac{f^{n-1}f'}{r^{n-1}}\right) 
=\left(\frac{\kappa}{n}r-\lambda r^{1-n}\right)f'. 
\label{eq:odeas} 
\eeq 
(By rescaling $f$ we can always set $\kappa=1$, but we choose 
to leave $\kappa$ arbitrary.) In the case $n=1$ it is straightforward to 
integrate (\ref{eq:ode}) to obtain the general solution 
$$ 
f= \frac{1}{\rho}
\log\left( \frac{\rho-\lambda+\kappa r}{\rho+\lambda-\kappa r}\right)+\alpha, 
$$ 
so that the potential is 
$$ 
g=\rho^{-1}(r-r_-) 
\log(r-r_-)+ 
\rho^{-1}(r_+-r)
\log(r_+-r)+\alpha r +\beta, 
\quad r_{\pm}=(\lambda\pm\rho)/\kappa,     
$$ 
for arbitrary constants $\rho>0$, $\alpha$, $\beta$.  
Up to shifting by an affine function, 
the solution for $n=1$ is of the form (\ref{eq:polys}), 
with $g(r)$ defined on a single interval 
$(r_-,r_+)$ 
in ${\mathbb R}$, which 
leads to a toric metric on $S^2$. 
Clearly for any $n$ there is the trivial solution 
$f=constant$, 
for which the inverse Hessian ${\bf G}^{-1}$  
in (\ref{eq:radmat}) becomes infinite. 
For $n>1$ we 
are unable to integrate  (\ref{eq:ode}) explicitly, 
and must resort to asymptotic analysis around singular points, 
before considering the solutions of the initial value problem. 

\noindent 
{\bf Asymptotics at $r=0$:} There are several different types of 
behaviour near the origin. 
\begin{itemize} 
\item For any $\lambda$, 
(\ref{eq:ode}) admits the expansion 
\beq 
f\sim f_0+ar^n -\frac{\lambda a^2n^2}{(n+1)}r^{n+1}+O(r^{n+2}), 
\label{const1} 
\eeq 
with arbitary constants $f_0\neq 0$ and $a$. 
\item For $\lambda\neq 0$ there is an alternative expansion 
\beq 
f\sim f_0+\frac{1}{\lambda(n-1)}r^{n-1}+O(r^{n+1}) 
\label{const2} 
\eeq 
with $f_0\neq 0$, for $n\neq 2$.  When $n=2$ we have instead
$$
f\sim f_0+\lambda^{-1}r-\frac{1}{4\lambda^2f_0}r^2+O(r^3).
$$
   
\item For $\lambda\neq 0$ and $n\neq 2$ there is 
another local expansion in  
(\ref{eq:ode}) that is regular at $r=0$, with 
leading order behaviour 
\beq 
f\sim-\frac{n(n-2)}{\lambda (n-1)}\, r^{n-1}+O(r^{n+1}). 
\label{eq:orig} 
\eeq 
for $n\neq 2$. 
\item In the case $\lambda=0$, there is an exact solution that is  singular at the 
origin, namely 
$$ 
f= \frac{2n^2}{\kappa}r^{-1}. 
$$ 
\end{itemize} 

\noindent 
{\bf Asymptotics at infinity:} 
The equation (\ref{eq:ode}) 
admits two types of asymptotic expansion at infinity: 
\begin{itemize} 
\item In the first, $f$ is asymptotic to a nonzero constant, i.e. 
\beq 
f\sim f_\infty +\frac{n(n+1)}{\kappa}r^{-1}+\ldots, 
\label{eq:asym} 
\eeq 
for $f_\infty\neq 0$. 
\item In the second case, $f$ tends to zero:  
\beq 
f\sim \frac{2n^2}{\kappa}r^{-1}, \qquad r\to\infty. 
\label{infzero} 
\eeq 
\end{itemize} 

\noindent
{\bf Movable singularities at $r=r_0\neq 0$:} The ODE 
(\ref{eq:ode}) is outside the Painlev\a'{e} class 
of second order equations whose general solution 
has no movable singularities 
other than poles, as described in \cite{ince}, 
since its solutions admit algebraic branching at arbitrary points 
$r_0\neq 0$ in the complex $r$ plane. The leading order 
behaviour in the neighbourhood of such a branch 
point is 
\beq 
f\sim c_0(r-r_0)^{1/n}, 
\label{eq:branch} 
\eeq 
where the constant $c_0$ is arbitrary. Since (\ref{eq:ode}) 
is second order, we expect that 
(\ref{eq:branch}) should be the first term in 
an expansion in powers of $(r-r_0)^{1/n}$ providing a  
local representation of the general solution, since it 
contains the two arbitrary constants $r_0$ and $c_0$. 
If it has only 
algebraic branching around movable singularities, an  
ODE can possess  
the weak Painlev\a'{e} property 
as defined in \cite{weak}, and there are many examples of 
integrable ordinary  and partial  
differential equations which have this property (see for instance 
\cite{abenda, hodo, hone}).  
However, the ODE (\ref{eq:ode}) also admits movable 
logarithmic branch points in its solutions, with 
\beq 
f\sim \left( \frac{\lambda}{r_0^{n-1}}-\frac{\kappa r_0}{n}\right)^{-1}\log (r-r_0) 
\label{eq:logb} 
\eeq 
in the neighbourhood of $r=r_0$. Logarithmic branching such 
as (\ref{eq:logb}) is taken as a 
strong indicator of non-integrability in differential equations 
\cite{ars}.  

We would expect that generically the solutions of the ODE should 
have infinitely many branch points in the complex plane. 
If a solution has a branch point on the real axis, then  
for even $n$ it cannot be real-valued for all 
real $r$. In that case, 
for $\lambda\neq 0$ it would make sense to consider 
a real-valued branch defined on the interval $[0,r^*)$, 
where $r^*>0$ is the position of the first positive real 
branch point. Such a solution might have the behaviour 
(\ref{eq:orig}) or tend to a nonzero constant at the origin, 
and as $r\to r^*$ from 
below would look like 
\beq 
f\sim c (r^*-r)^{1/n} \label{brabelow}  
\eeq 
or 
\beq 
f\sim \hat{c} \log (r^{*}-r) \label{logbelow}   
\eeq 
for real constants $c$, $\hat{c}$.       
In that case the Hessian matrix ${\bf G}$ would be entirely 
determined by the function $f$ and its first derivative, 
giving a metric on the ball of radius $r^*$ in  
${\mathbb R}^n$, with a singularity at $r=r^*$ 
where $f'$ blows up. Exactly how to extend this 
to a metric on a $2n$-dimensional (symplectic) manifold 
is not clear to us. 

In fact starting from 
initial data specified near the origin at $r=\epsilon>0$, 
the ODE (\ref{eq:ode}) is regular. It is not possible to 
take initial data at the origin since 
$r=0$ is a singular point of (\ref{eq:ode}). Clearly the right 
hand side of (\ref{eq:ode}) is also singular at $f=0$.  
It is easy to see 
that given $f'(\epsilon)\neq 0$, 
the  solution cannot have stationary points for $r>\epsilon$, since 
integrating (\ref{eq:ode}) with initial data $f(\epsilon)\neq 0$, 
$f'(\epsilon)=0$ gives only the constant solution $f(r)=f(\epsilon)$. It is 
then trivial to prove the following: 
\vspace{.1in} 

\noindent 
{\bf Lemma:} {\it Suppose that 
initial data is specified for the ODE (\ref{eq:ode}) 
at a point $r=\epsilon>0$, with $f(\epsilon)$, 
$f'(\epsilon)$ both nonzero. 
Then the solution $f(r)$ for $r>\epsilon$ remains positive (negative) 
for $f(\epsilon)>0$ ($f(\epsilon)<0$),  and for $f'(\epsilon)>0$ or 
$f'(\epsilon)<0$ it is  
monotone (increasing or decreasing, respectively),  
as long as it exists.}  
\vspace{.1in} 

We wish to consider solutions of the ODE which are  
free of branch points on the whole positive 
real axis, or on a finite 
interval $[0,r^*)$ (with  branching at $r=r^*$). In that case   
the components of the Hessian matrix ${\bf G}$ determine a metric 
on ${\mathbb R}^n$, or on the ball of radius $r^*$. The eigenvalues 
of ${\bf G}$ in (\ref{eq:radmat}) 
are $f/r$ (repeated $n-1$ times) and $f'$. 
Thus if we also require 
a Riemannian metric given by a positive definite Hessian, 
then we 
should consider only solutions of (\ref{eq:ode}) defined 
with both initial data positive, i.e. $f(\epsilon)>0$, 
$f'(\epsilon)>0$. The signs of $f$ and $f'$ remain constant  
and  the Lemma also holds when  the solution 
is analytically continued to $r<\epsilon$, provided it exists.  
However, a priori we have no guarantee 
that when the solution is continued back  towards $r=0$ (solving 
the initial value problem in reverse) it will not 
reach a singularity (branch point) at some point 
$r=\epsilon '$ with $0<\epsilon '<\epsilon$. 
Hence we are led to consider a nonlinear connection 
problem for (\ref{eq:ode}), 
to determine what asymptotic behaviours at 
$r=0$ are compatible with a solution defined on $[0,r^*)$ 
with specified branching at $r=r^*$, or a solution on the
whole positive real axis 
with specified 
asymptotics at infinity. 

For completeness we will also consider the cases when either 
$f$, $f'$ or both are negative, and henceforth we assume      
$\kappa>0$, $\lambda\neq 0$. Taking the four different 
combinations of initial data in turn, we obtain the   
following result: 
\vspace{.05in} 

\noindent 
{\bf Theorem 1:} {\it Suppose non-zero initial data is specified 
for the ODE (\ref{eq:ode}) at $r=\epsilon>0$. 
There are four possibilities: 

\noindent 
(i) $f(\epsilon)>0$, 
$f'(\epsilon)>0$. The solution reaches a logarithmic 
branch point with asymptotic behaviour (\ref{logbelow}) 
(with $\hat{c}<0$) 
at some finite point $r=r^*>\epsilon$. 

\noindent 
(ii) $f(\epsilon)<0$, 
$f'(\epsilon)>0$. The solution reaches an algebraic
branch point with asymptotic behaviour (\ref{brabelow}) 
(with $c<0$)   
at some finite point $r=r^*>\epsilon$. 

\noindent 
(iii) $f(\epsilon)>0$,
$f'(\epsilon)<0$. Either the solution reaches an algebraic  
branch point (\ref{brabelow}) (with $c>0$)  
at some finite point $r=r^*>\epsilon$, or $f$ is defined for all  
$r\geq\epsilon$ and asymptotes to a non-negative constant 
with the behaviour (\ref{eq:asym}) or (\ref{infzero}) 
as $r\to\infty$. 

\noindent 
(iv) $f(\epsilon)<0$,
$f'(\epsilon)<0$. Either the solution reaches a logarithmic
branch point (\ref{logbelow}) (with $\hat{c}>0$) 
at some finite point $r=r^*>\epsilon$, or $f$ is defined 
for all $r\geq\epsilon$ and asymptotes to a negative constant 
with the behaviour (\ref{eq:asym})  
as $r\to\infty$.    
} 
 
\noindent 
{\bf Proof.} (i) The solution cannot be defined 
for all $r\geq\epsilon$, since the only possible asymptotics 
(\ref{eq:asym}) or (\ref{infzero}) at infinity necessitate 
$f'<0$ for sufficiently large $r$ (with our assumption of 
positive $\kappa$), which contradicts the Lemma. 
An algebraic branch point (\ref{brabelow}) 
would mean $f\to 0$ as $r\to r^*$, contradicting the fact 
that $f$ is monotone increasing by the Lemma. Hence 
the logarithmic branching (\ref{logbelow}) at some finite 
$r=r^*$ is the only possibility. (ii) Again by the Lemma, 
$f$ must remain negative and monotone increasing. This 
rules out the asymptotic behaviours 
(\ref{eq:asym}) or (\ref{infzero}) 
at infinity for which $f'<0$, and also a logarithmic branch 
point which would require $f\to -\infty$. (iii) and (iv) are proved 
similarly. 
\vspace{.1in} 
         
Having considered the standard initial value problem for 
(\ref{eq:ode}) at $r=\epsilon$, we can now solve it in reverse, 
and consider continuing the solutions of types (i)-(iv) 
backwards towards the origin. (This is equivalent to solving 
the ODE for $h(r)=f(-r)$ with initial data at $r=-\epsilon$.) 
The solution to the nonlinear 
connection problem may be summarized thus:       
\vspace{.1in} 

\noindent 
{\bf Theorem 2:} {\it Solving the initial value problem 
for the ODE (\ref{eq:ode}) in the reverse direction, with the 
four combinations of initial 
data as in Theorem 1, leads to the following possibilities for 
the solution $f(r)$ with $r<\epsilon$: 

\noindent 
(i) Either the solution  reaches an algebraic branch point 
at $r=r_0$ with (\ref{eq:branch}) ($c_0>0$) 
for $0<r_0<\epsilon$, 
or it continues back to the origin with one of the 
asymptotic behaviours (\ref{const1}) (with 
$f_0>0$, $a>0$), 
(\ref{const2}) (with $f_0>0$, for $\lambda>0$ only) or 
(\ref{eq:orig}) (for $\lambda<0$, $n\neq 2$ only).          

\noindent
(ii) Either the solution  reaches a logarithmic branch point
at $r=r_0$ with (\ref{eq:logb})  
for $0<r_0<\epsilon$,
or it continues back to the origin with one of the
asymptotic behaviours (\ref{const1}) (with 
$f_0<0$, $a>0$),
or (\ref{const2}) (with $f_0<0$, for $\lambda>0$ only). 

\noindent
(iii) Either the solution  reaches a logarithmic branch point
at $r=r_0$ with (\ref{eq:logb}) for $0<r_0<\epsilon$,
or it continues back to the origin with one of the
asymptotic behaviours (\ref{const1}) (with 
$f_0>0$, $a<0$), 
or (\ref{const2}) (with $f_0<0$, for $\lambda<0$ only).

\noindent
(iv) Either the solution 
reaches an algebraic branch point
at $r=r_0$ with (\ref{eq:branch}) ($c_0<0$) 
for $0<r_0<\epsilon$,
or it continues back to the origin with one of the
asymptotic behaviours (\ref{const1}) (with 
$f_0<0$, $a<0$),
(\ref{const2}) (with $f_0<0$, for $\lambda<0$ only) or
(\ref{eq:orig}) (for $\lambda>0$, $n\neq 2$ only). 
} 

\noindent
{\bf Proof.} The proof is straightforward, by a simple 
enumeration of the various possibilities 
(\ref{const1}), (\ref{const2}) or (\ref{eq:orig}) at $r=0$ 
allowed by the Lemma, as well as the different branching 
behaviours that can occur.

\section{Conclusions} 

For any dimension $n>1$, new solutions of Abreu's equation 
(\ref{eq:abreueq}) have been found 
by making a similarity reduction to 
$O(n)$ invariant solutions. The similarity solutions 
are determined in terms of 
a single function $f(r)$ depending only on the radial 
coordinate $r$, which satisfies a second order 
ODE (\ref{eq:ode}) including an arbitrary parameter 
$\lambda$. Due to a singularity in the ODE at the origin, the 
initial value problem is not defined at $r=0$. Nevertheless 
we have shown that for some suitable subset of the 
possible initial data 
specified at $r=\epsilon>0$ the solution  
may have a continuation in both the forward ($r>\epsilon$) 
and backward ($r<\epsilon$) directions to 
define a solution either on the whole positive real axis, 
or on an interval $[0,r^*)$ with a branch point at 
$r=r^*$. With the restriction to  Riemannian        
metrics defined by a positive definite Hessian matrix 
${\bf G}$ in (\ref{eq:radmat}), only the 
latter case is relevant (corresponding to 
case (i) in Theorems 1 and 2 above). However, whether such 
solutions could have further geometric significance, 
by extension to a suitable metric on a symplectic 
manifold of dimension $2n$, is uncertain. 

Performing numerical integrations of (\ref{eq:ode}) 
with $n=3$ and $\kappa=\lambda=1$, by different choices 
of initial data  
we have obtained solutions displaying some of the 
asymptotic behaviours included in cases (i)-(iv) 
of the Theorems.  
A more 
detailed understanding of such solutions would 
necessitate estimates on the initial data to 
ensure that branch points do not appear as  
$r\to 0$. Determining conditions for the existence of 
points of inflection would be 
an essential step in such an understanding.         
It would also be desirable to have some bounds 
on $r^*$, the position of the branch point, in terms 
of $f(\epsilon)$ and $f'(\epsilon)$, but this is beyond the 
scope of the present work.  

The various Painlev\'{e} tests are good heuristic tools 
for analysing ordinary and partial differential 
equations, and have been used to isolate new 
integrable systems (see Chapter 7 of \cite{abcla} for 
a review). These tests are based on the Painlev\'{e} 
property, that the general solution of an equation should 
have only poles as movable singularities. However, as 
originally emphasized by the authors of \cite{ars}, 
the Painlev\'{e} property is extremely sensitive to 
changes of variables. For example, transformations 
of hodograph type can change equations with movable 
algebraic branch points into equations that have only 
poles \cite{hodo, hone}. The ODE (\ref{eq:ode}) has 
not only algebraic branching but also movable 
logarithmic branch points, and so certainly fails 
the standard Painlev\'{e} test or its weak 
extension \cite{weak}. This is a strong indication 
that Abreu's equation (\ref{eq:abreueq}) is not 
integrable. 

It has been observed \cite{smm} that when an 
equation is written in potential form it can 
have a single (or finitely many) logarithmic terms 
in a local expansion around singular points, and still be integrable.    
However, the asymptotic expansion for (\ref{eq:ode}) 
with (\ref{eq:logb}) as the first term seems to require 
infinitely many logarithms. In the case of the  
KP equation (\ref{kp}), the potential $w$ in (\ref{kplax}) 
is given by 
the logarithmic derivative of the tau-function: 
$$ 
w=-\frac{\partial}{\partial x}\log \tau. 
$$ 
The tau-function for the rational solutions (\ref{kprat}) 
takes the form  
\beq 
\tau=\prod_{l=1}^d (x-x_l(\underline{t})), 
\label{polytau} 
\eeq   
and this is a polynomial both in $x$ and the times $\underline{t}$. 
For the more general algebro-geometric solutions of KP 
\cite{krich},  the  
tau-function is a theta-function of an arbitrary 
Riemann surface. The rational and soliton solutions 
arise as degenerate limits of the theta-functions. 
It would be interesting to see if Abreu's equation 
would admit quasiperiodic generalizations of 
(\ref{absum}), 
by dropping the extremality condition (\ref{extrem}) on the  
the curvature, allowing $S$ in (\ref{eq:curv}) to 
be a more general (say, periodic) function.   
  
Other methods of testing for integrability might  
provide useful information about Abreu's equation. 
Since (\ref{eq:abreueq}) is not of evolution type, 
a promising method would be 
the symmetry approach of Shabat et al \cite{shabat},  
which has recently been extended \cite{mik} 
in order to deal with non-evolutionary equations.  
It would also be 
interesting to apply the methods of \cite{olver} 
to look for other sorts of group invariant solutions 
of Abreu's equation (\ref{eq:abreueq}).

\end{document}